\documentclass[thmsa]{article}
\usepackage{amsfonts}
\usepackage{amssymb}
\usepackage{amsmath}
\usepackage{amsmath}
\usepackage{graphicx}

\begin{document}

\begin{center}
{\Large Generalized Iterated Function System and Common Attractors of
Generalized Hutchinson Operators in Dislocated Metric Spaces}

Talat Nazir$^{1}$ and Sergei Silvestrov$^{2}$

$^{1}$Department of Mathematical Sciences, University of South Africa,

Florida 0003, South Africa

talatn@unisa.ac.za

$^{2}$Division of Mathematics and Physics, School of Education, Culture and
Communication,

M\"{a}lardalen University, Box 883, 72123 V\"{a}steras, Sweden

sergei.silvestrov@mdu.se{\large \medskip }

--------------------------------------------------------------------------------------------
\end{center}

\noindent {\footnotesize \textbf{Abstract. }}In this paper, we present the
generalized iterated function system for constructing of common fractals of
generalized contractive mappings in the setup of dislocated metric spaces.
The well-posedness of attractors based problems of rational contraction maps
in the framework of dislocated metric spaces is also established. Moreover,
the generalized collage theorem is also established in dislocated metric
spaces.


\noindent {\footnotesize \textbf{AMS Classification 2010:} 47H10, 54H25,
54E50}

\section{Introduction and preliminaries}

Metric fixed point theory serves as an essential tool for solving problems
arising in various branches of mathematical analysis for instance, split
feasibility problems, variational inequality problems, nonlinear
optimization problems, equilibrium problems, complementary problems,
selection and matching problems, and problems of proving an existence of
solution of integral and differential equations, adaptive control systems,
fractal image decoding, convergence of recurrent networks and many more. In
particular, it has deep roots in nonlinear functional analysis.

Hitzler and Seda \cite{Hitzler2000}\ introduced the notion of dislocated
metric space and proved a fixed point result as an interesting
generalization of the Banach contraction principle. Then several useful
results of fixed point theory were established, \cite{Abbas19, Aydi18,
Karapynar2013, Sumati2015} introduced the notion of contractions and related
fixed point results in dislocated metric spaces.

By using the concept of Hausdorff dislocated metric and generalized
contraction mappings, we establish the existence of attractors of Hutchinson
operators in the framework of Dislocated metric spaces. First we recall the
definition of dislocated metric space.

Thought this work, set of real numbers shall be represented by $%
\mathbb{R}
$, set of non-negative real numbers by $%
\mathbb{R}
_{+}$, set of $a$-tuples of real numbers by $%
\mathbb{R}
^{a}$ and set of natural numbers by $%
\mathbb{N}
$. First, we review some key concepts.\medskip

\noindent \textbf{Definition 1.1.} \cite{Hitzler2000}\ Let $X$ be a nonempty
set. A function $\delta :X\times X\rightarrow 
\mathbb{R}
_{+}$ is said to be a dislocated metric (or a metric -like) on $X$ if for
any $x,y,z\in X,$ the following condition hold:

\begin{enumerate}
\item $\delta (x,y)=0$ implies that $x=y;$

\item $\delta (x,y)=\delta (y,x);$

\item $\delta (x,z)\leq \delta (x,y)+\delta (y,z).$
\end{enumerate}

\noindent Then $\delta $ is called a dislocated metric and the pair $%
(X,\delta )$ is called a dislocated metric space.\medskip

\noindent \textbf{Example 1.2. \  \ }We take $X=\{a,b,c\} \subseteq 
\mathbb{R}
$ and consider the dislocated metric $\delta :X\times X\rightarrow 
\mathbb{R}
_{+}$ given by%
\begin{eqnarray*}
\delta (a,a) &=&0,\text{ \  \  \ }\delta (b,b)=1,\text{ \  \  \ }\delta (c,c)=%
\frac{2}{3}, \\
\delta (a,b) &=&\delta (b,a)=\frac{9}{10},\text{ \  \  \ }\delta (b,c)=\delta
(c,b)=\frac{4}{5}, \\
\delta (a,c) &=&\delta (c,a)=\frac{7}{10}.
\end{eqnarray*}%
Since $\delta (b,b)\neq 0,$ $\delta $ is not a metric and since $\delta
(b,b)\geq \delta (a,b),$ $\delta $ is not a partial metric defined in \cite%
{Aydi18}.\medskip

\noindent \textbf{Example 1.3. \  \ }Let $X=%
\mathbb{R}
_{+}$ and $a,b\in 
\mathbb{R}
_{+}.$ Consider $\delta :X\times X\rightarrow 
\mathbb{R}
_{+}$ given by%
\begin{equation*}
\delta (x,y)=a\left \vert x-y\right \vert +b\max \{x,y\}.
\end{equation*}%
If we take $a=1$ and $b=0$, then $\delta $ is a metric on $X$.

\noindent If we take $a=0$ and $b=1$, then $\delta $ is a partial metric on $%
X$.

\noindent If we take $a=2$ and $b=4$, then $\delta $ is a dislocated metric
on $X$ which is neither metric nor partial metric on $X.$\medskip 

\noindent \textbf{Definition 1.4 }(Open Ball). Let $(X,\delta )$ be a
dislocated metric space and $\varepsilon \in 
\mathbb{R}
_{+}.$ We define the open ball as follows:%
\begin{equation*}
B_{\varepsilon }(x)=\{y\in X:|\delta (x,y)-\delta (x,x)|<\varepsilon \}.
\end{equation*}%
The topology $\tau _{\delta }$ on $(X,\delta )$ as follows:%
\begin{equation*}
\tau _{\delta }=\{U\subseteq X:\text{ }\forall \text{ }u\in U\text{ }\exists 
\text{ }\varepsilon \in 
\mathbb{R}
_{+}\text{ such that }B_{\varepsilon }(x)\sqsubseteq U\}.
\end{equation*}%
A sequence $\{u_{n}\}$ in$\ (X,\delta )$ is said to convergence to $u\in X$
if and only if $\lim \limits_{n\rightarrow \infty }\delta (u_{n},u)=\delta
(u,u).$

A limit of every convergence sequence in dislocated metric space $(X,\delta )
$ is unique \cite{Zeyada06}.\medskip 

\noindent \textbf{Definition\ 1.5.} \cite{Hitzler2000}\ Let $(X,\delta )$ be
a dislocated metric space.

\begin{enumerate}
\item A sequence $\{x_{n}\}$ in $X$ is said to be Cauchy sequence if $%
\underset{n,m\rightarrow \infty }{\lim }\delta (x_{n},x_{m})$ exists and is
finite.

\item $(X,\delta )$ is said to be complete if every Cauchy sequence $%
\{x_{n}\}$ in $X$ converges with respect to $\tau _{\delta }$ to a point $%
x\in X$ such that $\underset{n\rightarrow \infty }{\lim }\delta
(x_{n},x)=\delta (x,x)=\underset{n\rightarrow \infty }{\lim }\delta
(x_{n},x_{m}).$
\end{enumerate}

\noindent \textbf{Definition\ 1.6.} \cite{Zeyada06}\ Let $(X_{1},\delta _{1})
$ and $(X_{2},\delta _{2})$ be two dislocated metric spaces and let\ $%
\mathfrak{f}:X_{1}\rightarrow X_{2}$ be a function. $\mathfrak{f}$ is said
to be continuous if for each sequence $\{u_{n}\}$ which converges to $u_{0}$
in $X_{1}$, then sequence $\{ \mathfrak{f}\left( u_{n}\right) \}$ converges
to $\mathfrak{f}\left( u_{0}\right) $ in $X_{2}$.\medskip 

A subset $Y$ in\ dislocated metric space $(X,\delta )$ is said to be bounded
if and only if the set $\{ \delta (x,y):x,y\in Y\}$ is bounded above.\medskip 

\noindent \textbf{Definition\ 1.7. }\cite{Aydi18} Let $\overline{C}$ be a
closure of $C$ with respect to dislocated metric $\delta .$ Then%
\begin{equation*}
c\in \overline{C}\Longleftrightarrow B_{\delta }(c,\varepsilon )\cap C\neq
\emptyset \text{ \  \ for all }\varepsilon >0.
\end{equation*}%
A set $C$ in dislocated metric space is closed if and only if$\  \overline{C}%
=C.$\medskip 

\noindent \textbf{Definition\ 1.8.} \ $\ (X,\delta )$ be a dislocated metric
space. A subset $K$ of $X$ is said to be compact if and only if every open
cover of $K$ (by open sets in $M$) has a finite subcover. If $M$ itself has
this property, then we say that $M$ is a compact dislocated metric
space.\medskip 

\noindent \textbf{Theorem\ 1.9.} \ $\ $Let $(X,\delta )$ be a dislocated
metric space, and let $K$ be a compact subset of $X$. Then $K$ is a closed
subset of $X$, and $K$ is bounded.

\noindent \textbf{Proof}.\  \  \ Let $K$ be a compact subset of a dislocated
metric space $(X,\delta ).$ To show that $K$ is closed, we show that the
complement, $O=X\backslash K$, is open. Let $z\in O$. We need to find $%
\varepsilon >0$ such that $B\varepsilon (z)\subseteq O$. Now, for any $x\in K
$, let $\varepsilon _{x}=d(x,z)$. Since $z\notin K$, $\varepsilon _{x}>0$.
The collection of open sets $\{B\varepsilon _{x}(x):x\in K\}$ is an open
cover of $K$ (since any $x\in K$ is covered by $B\varepsilon _{x}(x)$).
Since $K$ is compact, there is a finite subcover of this cover, that is,
there is a finite set $x_{1},x_{2},...,x_{n}$ such that the corresponding
open balls already cover $K$. Let $\varepsilon =\frac{1}{2}\min
\{ \varepsilon _{x_{1}},\varepsilon _{x_{2}},...,\varepsilon _{x_{n}}\}$. The
claim is that $B_{\varepsilon }(z)\subseteq O$. To show this, let $y\in
B_{\varepsilon }(z)$. We want to show $y\in O$, that is, $y\notin K$.
Consider $x_{i}$, where $1\leq i\leq n$. By the triangle inequality, $\delta
(x_{i},y)+\delta (y,z)\geq \delta (x_{i},z)=\varepsilon _{x_{i}}\geq
2\varepsilon $. So, $\delta (x_{i},y)\geq 2\varepsilon -\delta
(y,z)>2\varepsilon -\varepsilon =\varepsilon $. (The last inequality follows
because $\delta (y,z)<\varepsilon $). Then, since $\delta
(x_{i},y)>\varepsilon $, $y$ is not in the open ball of radius $\varepsilon
_{x_{i}}$ about $x_{i}$. Since the open balls $B\varepsilon _{x_{i}}(x_{i})$
cover $K$, we have that $y\notin K$.

Now, to prove that $K$ is bounded, let $x$ be an element of $K$, and
consider the collection of open balls of integral radius, $%
\{B_{i}(x):i=1,2,...\}$. Since every element of $K$ has some finite distance
from $x$, this collection is an open cover of $K$. Since $K$ is compact, it
has a finite subcover $\{B_{i_{1}}(x),B_{i_{2}}(x),...,B_{i_{n}}(x)\}$,
where we can assume $i_{1}<i_{2}<$\textperiodcentered \textperiodcentered
\textperiodcentered $<i_{n}$. But since $B_{i_{1}}(x)\subseteq
B_{i_{2}}(x)\subseteq $\textperiodcentered \textperiodcentered
\textperiodcentered $\subseteq B_{i_{n}}(x)$, this means that $B_{i_{n}}$ by
itself already covers $K$. Then, for $y,z\in K$, $y$ and $z$ are in $%
B_{x_{n}}(x)$, and $\delta (y,z)\leq \delta (y,x)+\delta (x,z)\leq
i_{n}+i_{n}$. It follows that $2i_{n}$ is an upper bound for $\{ \delta
(y,z):y,z\in K\}$. Thus, $K$ is bounded.\medskip 

\noindent \textbf{Definition 1.10.} \  \ A dislocated metric space $(X,\delta
)$ metric space is sequentially compact if every sequence has a convergent
subsequence.\medskip 

\noindent \textbf{Theorem 1.11}. \  \ A dislocated metric space $(X,\delta )$
is compact if and only if it is sequentially compact.

\noindent \textbf{Proof.} \  \ Suppose that $X$ is compact. Let $\{F_{n}\}$
be a decreasing sequence of closed nonempty subsets of $X$, and let $%
G_{n}=F_{n}^{c}$.

\noindent If $\cup _{n=1}^{\infty }G_{n}=X$, then $\{G_{n}:n\in 
\mathbb{N}
\}$ is an open cover of $X$, so it has a finite subcover $%
\{G_{n_{k}}:k=1,2,...,l\}$ since $X$ is compact. Let $N=\max
\{n_{k}:k=1,2,...l\}$.

\noindent Then $\cup _{n=1}^{N}G_{n}=X$, so $F_{N}=\cap
_{n=1}^{N}F_{n}=\left( \cup _{n=1}^{N}G_{n}\right) ^{c}=\emptyset $,
contrary to our assumption that every $F_{n}$ is nonempty.

\noindent It follows that $\cup _{n=1}^{\infty }G_{n}\neq X$ and then $\cap
_{n=1}^{\infty }F_{n}=\left( \cup _{n=1}^{\infty }G_{n}\right) ^{c}\neq
\emptyset ,$ meaning that $X$ has the finite intersection property for
closed sets, so $X$ is sequentially compact.

\noindent Conversely, suppose that $X$ is sequentially compact. Let $%
\{G\alpha \subset X:\alpha \in I\}$ be an open cover of $X$. Then, there
exists $\eta >0$ such that every ball $B_{\eta }(x)$ is contained in some $%
G\alpha $. Since $X$ is sequentially compact, it is totally bounded, so
there exists a finite collection of balls of radius $\eta ,$%
\begin{equation*}
\{B\eta (x_{i}):i=1,2,...,n\}
\end{equation*}%
that covers $X$. Choose $\eta _{i}\in I$ such that $B\eta (x_{i})\subset
G\alpha _{i}$. Then $\{G\alpha _{i}:i=1,2,...,n\}$ is a finite subcover of $X
$, so $X$ is compact.\medskip 

\noindent \textbf{Theorem 1.12}. \  \ Let $f$ be a continuous selfmap on
compact set $X$ in dislocated metric space $(X,\delta $) into itself. Then
the range $f(X)$ of $f$ is also compact.

\noindent \textbf{Proof}. \  \ We are to show that, for any sequence $%
\{a_{n}\}$ in $X,$ the sequence $\{fa_{n}\}$ has a convergent subsequence
with limit $\{fa_{0}\}$ for some $a_{0}$ in $X.$

\noindent Since the sequence $\{a_{n}\}$ in $X,$ we have a subsequence $%
\{a_{n_{k}}\}$ of $\{a_{n}\}$ such that $\{a_{n_{k}}\}$ converges to some $%
a_{0}$ in $X.$

\noindent By continuity of $f,$ we obtain that $\{fa_{n_{k}}\}$ converges to 
$\{fa_{0}\}.$

\noindent As the sequence $\{fa_{n}\}$ is in $X$ has a convergent
subsequence $\{fa_{n_{k}}\}$ with limit $\{fa_{0}\}$ for some $a_{0}$ in $X.$
Consequently, the range $f(X)$ of $f$ is also compact.\medskip 

In the dislocated metric space $\left( X,\delta \right) ,$ we define the
following sets:%
\begin{eqnarray*}
\mathcal{B}\left( X\right)  &=&\{Y:\text{ }Y\text{ is nonempty closed and
bounded subset of }X\}, \\
\mathcal{C}\left( X\right)  &=&\{Y:\text{ }Y\text{ is nonempty compact
subset of }X\}.
\end{eqnarray*}%
\noindent For $R,S\in \mathcal{B}(X)$ and $x\in X$, we define%
\begin{eqnarray*}
\delta (x,R) &=&\inf \{ \delta (x,r),r\in R\} \text{ and} \\
\sigma _{\delta }(R,S) &=&\sup \{ \delta (r,S):r\in R\}.
\end{eqnarray*}%
\noindent \textbf{Definition 1.13. }\cite{Aydi18}\textbf{\ }Let $(X,\delta )$
be a dislocated metric space. For $R,S\in \mathcal{B}(X),$ we define
dislocated Hausdroff metric $H_{\delta }$ by%
\begin{equation*}
H_{\delta }(R,S)=\max \{ \sigma _{\delta }(R,S),\sigma _{\delta }(S,R)\}.
\end{equation*}%
The pair $(X,H_{\delta })$ is called dislocated Hausdroff metric
space.\medskip 

\noindent \textbf{Theorem 1.14. }\cite{Aydi18} Let $(X,\delta )$\textbf{\ }%
be a dislocated metric space. For all $U,V,W\in \mathcal{B}\left( X\right) .$
Then

\begin{description}
\item[(H1)] $H_{\delta }(R,S)=\sigma _{\delta }(R,R)=\sup \{ \sigma
(r,R):r\in R\};$

\item[(H2)] $H_{\delta }(R,S)=H_{\delta }(S,R);$

\item[(H3)] $H_{\delta }(R,S)=0\Longrightarrow S=R;$

\item[(H3)] $H_{\delta }(R,T)\leq H_{\delta }(R,S)+H_{\delta }(S,T).$
\end{description}

If $(X,\delta )$ is a complete dislocated metric space, then $\left( 
\mathcal{C}(X),H_{\delta }\right) $ is also complete dislocated Hausdroff
metric space.\medskip 

\noindent \textbf{Lemma 1.15.} \  \ Let $(X,\delta )$ be a dislocated metric
space. For all $R,S,U,V\in \mathcal{B}(X)$, the following hold:

\begin{description}
\item[(i)] If $S\subseteq U,$ then $\sup \limits_{r\in R}\delta (r,U)\leq
\sup \limits_{r\in R}\delta (r,S).$

\item[(ii)] $\sup \limits_{x\in R\cup S}\delta (x,U)=\max \{ \sup
\limits_{r\in R}\delta (r,U),\sup \limits_{s\in S}\delta (s,U)\}.$

\item[(iii)] $H_{\delta }(R\cup S,U\cup V)\leq \max \{H_{\delta
}(R,U),H_{\delta }(S,V)\}.$
\end{description}

\noindent \textbf{Proof}\textit{.} \  \ To prove (i):\ Since $S\subseteq U,$
for all $r\in R,$ we have%
\begin{eqnarray*}
\delta (r,U) &=&\inf \{ \delta (r,u):u\in U\} \\
&\leq &\inf \{ \delta (r,s):s\in S\}=\delta \left( r,S\right) ,
\end{eqnarray*}%
which implies that%
\begin{equation*}
\sup \limits_{r\in R}\delta (r,U)\leq \sup \limits_{r\in R}\delta (r,S).
\end{equation*}%
To prove (ii):%
\begin{eqnarray*}
\sup_{x\in R\cup S}\delta \left( x,U\right)  &=&\sup \{ \delta \left(
x,U\right) :x\in R\cup S\} \\
&=&\max \{ \sup \{ \delta \left( x,U\right) :x\in R\},\sup \{ \delta \left(
x,U\right) :x\in S\} \\
&=&\max \{ \sup_{r\in R}\delta \left( r,U\right) ,\sup_{s\in S}\delta \left(
s,U\right) \}.
\end{eqnarray*}%
\noindent To prove (iii): Note that%
\begin{eqnarray*}
&&\sup \limits_{x\in R\cup S}\delta (x,U\cup V) \\
&\leq &\max \{ \sup \limits_{r\in R}\delta (r,U\cup V),\sup \limits_{s\in
S}\delta (s,U\cup V)\} \text{ (by using (ii))} \\
&\leq &\max \{ \sup \limits_{r\in R}\delta (r,U),\sup \limits_{s\in S}\delta
(s,V)\} \text{ (by using (i))} \\
&\leq &\max \left \{ \max \{ \sup \limits_{r\in R}\delta
(r,U),\sup \limits_{u\in U}\delta (u,R)\},\max \{ \sup \limits_{s\in S}\delta
(s,V),\sup \limits_{v\in V}\delta (v,S)\} \right \}  \\
&=&\max \left \{ H_{\delta }\left( R,U\right) ,H_{\delta }\left( S,V\right)
\right \} .
\end{eqnarray*}%
In the similar way, we obtain that%
\begin{equation*}
\sup \limits_{y\in U\cup V}\delta (y,R\cup S)\leq \max \left \{ H_{\delta
}\left( R,U\right) ,H_{\delta }\left( S,V\right) \right \} .
\end{equation*}%
Hence it follows that%
\begin{eqnarray*}
H_{\delta }(R\cup S,U\cup V) &=&\max \{ \sup \limits_{x\in R\cup S}\delta
(x,U\cup V),\sup \limits_{y\in U\cup V}\delta (y,R\cup S)\} \\
&\leq &\max \left \{ H_{\delta }\left( R,U\right) ,H_{\delta }\left(
S,V\right) \right \} .
\end{eqnarray*}%
\noindent \textbf{Definition 1.16. \  \ }Let $(X,\delta )$ be a dislocated
metric space and $f,g:X\rightarrow X$ be mappings. A pair of mappings $%
\left( f,g\right) $ is called generalized contraction if%
\begin{equation*}
\delta \left( fx,gy\right) \leq \alpha \delta \left( x,y\right) 
\end{equation*}%
for all $x,y\in X,$ where $0\leq \alpha <1$.\medskip 

\noindent \textbf{Theorem 1.17. \  \ }Let $(X,\delta )$ be a dislocated
metric space and $f,g:X\rightarrow X$ be two continuous mappings. If the
pair of mappings $\left( f,g\right) $ is generalized contraction with $0\leq
\alpha <1$. Then

\begin{description}
\item[(1)] the elements in $\mathcal{C}(X)$ are mapped to elements in $%
\mathcal{C}(X)$ under $f$ and $g.$

\item[(2)] if for any $U\in \mathcal{C}(X),$%
\begin{eqnarray*}
f(U) &=&\{f(u):u\in U\} \text{ and} \\
g(U) &=&\{g(u):u\in U\}.
\end{eqnarray*}%
Then for $f,g:\mathcal{C}(X)\rightarrow \mathcal{C}(X),$ the pair $\left(
f,g\right) $ is a generalized contraction map on $(\mathcal{C}(X),H_{\delta
})$.
\end{description}

\noindent \textbf{Proof}. To prove \textbf{(1)}: Since $f$ is continuous
mapping and the image of a compact subset under $f:X\rightarrow X$\ is
compact, that is,%
\begin{equation*}
U\in \mathcal{C}(X)\text{ implies }f(U)\in \mathcal{C}(X).
\end{equation*}%
Similarly, we have%
\begin{equation*}
U\in \mathcal{C}(X)\text{ implies }g(U)\in \mathcal{C}(X).
\end{equation*}%
To prove \textbf{(2)}: Let $A_{1},A_{2}\in \mathcal{C}(X)$. Since for $%
f,g:X\rightarrow X,$ the pair of mappings $\left( f,g\right) $ is a
generalized contraction, we obtain that%
\begin{equation*}
\delta \left( fu,gv\right) \leq \alpha \delta \left( u,v\right) \text{ for
all }u,v\in X,
\end{equation*}%
where $0\leq \alpha <1$.

\noindent Thus we have%
\begin{eqnarray*}
\delta \left( fu,g\left( V\right) \right) &=&\inf_{v\in V}\delta \left(
fu,gv\right) \\
&\leq &\inf_{v\in V}\alpha \delta \left( u,v\right) \\
&=&\alpha \delta \left( u,V\right) .
\end{eqnarray*}%
Also%
\begin{equation*}
\delta \left( fv,g\left( U\right) \right) =\inf_{u\in U}\delta \left(
fv,gu\right) \leq \inf_{u\in U}\alpha \delta \left( v,u\right) =\alpha
\delta \left( v,U\right) .
\end{equation*}%
Now%
\begin{eqnarray*}
H_{\delta }\left( f\left( U\right) ,g\left( V\right) \right) &=&\max \{ \sup
\limits_{u\in U}\delta (fu,f\left( V\right) ),\sup \limits_{v\in V}\delta
(fv,f\left( U\right) )\}, \\
&\leq &\max \{ \sup \limits_{u\in U}\alpha \delta (u,V),\sup \limits_{v\in
V}\alpha \delta (v,U)\}, \\
&=&\max \{ \alpha (\sup \limits_{u\in U}\delta (u,V)),\alpha (\sup
\limits_{v\in B}\delta (v,U))\}, \\
&=&\alpha \max \{ \sup \limits_{u\in U}\delta (u,V),\sup \limits_{v\in
V}\delta (v,U)\}, \\
&=&\alpha H_{\delta }\left( U,V\right) .
\end{eqnarray*}%
Consequently,%
\begin{equation*}
H_{\delta }\left( f\left( U\right) ,g\left( V\right) \right) \leq \alpha
H_{\delta }\left( U,V\right) .
\end{equation*}%
Hence, the pair $\left( f,g\right) $ is a generalized contraction map on $(%
\mathcal{CB}(X),H_{\delta })$. $\square $

\section{Generalized Iterated Function System}

In this section, we construct a fractal set of generalized iterated function
system, a certain finite collection of mappings defined in the setup of
dislocated metric space. We also define Hutchinson operator with the help of
a finite collection of generalized rational contraction mappings on a
dislocated metric space. We start with the following result.

\noindent \textbf{Proposition 2.1.} \  \ Let $(X,\delta )$ be a dislocated
metric space.\ Suppose that the mappings $f_{n},g_{n}:X\rightarrow X$ for $%
n=1,2,...,N\ $are satisfying%
\begin{equation*}
\delta \left( f_{n}x,g_{n}y\right) \leq \alpha _{n}\delta \left( x,y\right) 
\text{ for all }x,y\in X,
\end{equation*}%
where\ $0\leq \alpha _{n}<1$ for each $n\in \left \{ 1,2,...,N\right \} .$
Then the mappings $T,S:\mathcal{C}(X)\rightarrow \mathcal{C}(X)$ defined as%
\begin{eqnarray*}
T(U) &=&f_{1}(U)\cup f_{2}(U)\cup \cdot \cdot \cdot \cup f_{N}(U) \\
&=&\cup _{n=1}^{N}f_{n}(U),\text{ for each }U\in \mathcal{C}(X)
\end{eqnarray*}%
and%
\begin{eqnarray*}
S(U) &=&g_{1}(U)\cup g_{2}(U)\cup \cdot \cdot \cdot \cup g_{N}(U) \\
&=&\cup _{n=1}^{N}g_{n}(U),\text{ for each }U\in \mathcal{C}(X)
\end{eqnarray*}%
are also satisfy%
\begin{equation*}
H_{\delta }\left( TU,SV\right) \leq \alpha _{\ast }H_{\delta }\left(
U,V\right) \text{ for all }U,V\in \mathcal{C}(X),
\end{equation*}%
where $\alpha _{\ast }=\max \{ \alpha _{i}:i\in \{1,2,...,N\} \},$ that is,
the pair $\left( T,S\right) $ is a generalized contraction on $\mathcal{C}%
\left( X\right) $.

\noindent \textbf{Proof.}\  \  \ Need complete proof. $N=1$

Assume that its true for $N=k.$

Then prove it is true $N=k+1.$

We will prove the result for $N=2$. Let $f_{1},f_{2},g_{1},g_{2}:X%
\rightarrow X$ be two contractions. For $A_{1},A_{2}\in \mathcal{C}\left(
X\right) \ $and using Lemma (1.9) (iii), we have%
\begin{eqnarray*}
H_{\delta }(T\left( A_{1}\right) ,S(A_{2})) &=&H_{\delta }(f_{1}(A_{1})\cup
f_{2}(A_{1}),g_{1}(A_{2})\cup g_{2}(A_{2})) \\
&\leq &\max \{H_{\delta }(f_{1}(A_{1}),g_{1}(A_{2}))\},H_{\delta
}(f_{2}(A_{1}),f_{2}(A_{2}))\} \\
&\leq &\max \{ \alpha _{1}H_{\delta }(A_{1},A_{2}),\alpha _{2}H_{\delta
}(A_{1},A_{2})\} \\
&=&\alpha _{\ast }H_{\delta }(A_{1},A_{2}).\text{ \ }\square
\end{eqnarray*}

\noindent \textbf{Definition 2.2.} \  \ Let $(X,\delta )$ be a dislocated
metric space\ and $T,S:\mathcal{C}\left( X\right) \rightarrow \mathcal{C}%
\left( X\right) $ are two mappings. A pair of mappings $\left( T,S\right) $
is called a generalized rational contractive, if there $\alpha \in \lbrack
0,1)$ such that%
\begin{equation*}
H_{\delta }(T\left( U\right) ,S(V))\leq \alpha M_{T,S}(U,V),
\end{equation*}%
where%
\begin{eqnarray*}
M_{T,S}(U,V) &=&\max \{H_{\delta }(U,V),H_{\delta }(U,T\left( U\right)
),H_{\delta }(V,S\left( V\right) ), \\
&&\frac{H_{\delta }(U,S\left( V\right) )+H_{\delta }(V,T\left( U\right) )}{2}%
,\dfrac{H_{\delta }(V,S\left( V\right) )[1+H_{\delta }(V,T(U))]}{1+H_{\delta
}\left( U,V\right) }, \\
&&\frac{H_{\delta }(V,S\left( V\right) )[1+H_{\delta }(U,T\left( U\right) )]%
}{1+H_{\delta }(U,V)},\frac{H_{\delta }(V,T\left( U\right) )[1+H_{\delta
}(V,T\left( U\right) )]}{1+H_{\delta }(U,V)}, \\
&&\frac{H_{\delta }(V,T\left( U\right) )[1+H_{\delta }(U,T\left( U\right) )]%
}{1+H_{\delta }(U,V)}\}.
\end{eqnarray*}%
\noindent The above defined operator $(T,S)$ is also known as a generalized
rational contractive Hutchinson operator which is the extension of
Hutchinson operator given in \cite{Hutchinson81}. Moreover if $(T,S)$
defined in Proposition 2.1 is generalized contraction, then it is trivially
generalized rational contraction and so $(T,S)$ is generalized rational
contractive Hutchinson operator.\medskip

\noindent \textbf{Definition 2.3.} \  \ Let $X$ be a dislocated metric space.
If $f_{n},g_{n}:X\rightarrow X$, $n=1,2,...,N$ are generalized contraction
mappings, then $(X;f_{1},f_{2},...,f_{N};g_{1,}g_{2},...,g_{N})$ is called
generalized iterated function system (GIFS).\medskip

\noindent \textbf{Definition 2.4.} \  \ A nonempty closed and bounded subset $%
U$ of $X$ is called a common attractor of $T$ and $S\ $generated by GIFS if

\begin{itemize}
\item[(a)] $T(U)=S(U)=U$ and

\item[(b)] there exist an open subset $V$ of $X$ such that $U\subseteq V$
and $\lim \limits_{k\rightarrow \infty }T^{k}(B)=\lim \limits_{k\rightarrow
\infty }S^{k}(B)=U$ for any closed and bounded subset $B$ of $V$, where the
limit is applied with respect to the dislocated Hausdorff metric.
\end{itemize}

\noindent The largest open set $V$ satisfying (b) is called a basin of
attraction.

\section{Main Results}

By using the concept of Hausdorff dislocated metric and generalized rational
contractive Hutchinson operator, we establish the existence of unique common
attractors of such mappings satisfying more general contractive conditions
than those in \cite{Abbas19, Balaji2017, Karapynar2013} in the framework of
complete dislocated metric spaces. We start with the following Theorem.

\noindent \textbf{Theorem 3.1}. \  \ Let $(X,\delta )$ be a complete
dislocated metric space with 

\noindent $(X;f_{1},f_{2},...,f_{N};g_{1},g_{2},...,g_{N})$ be a given
generalized iterated function system.$\ $Suppose that the pair of self
mappings $(T,S)$ defined by%
\begin{eqnarray*}
T(B) &=&f_{1}(B)\cup f_{2}(B)\cup \cdot \cdot \cdot \cup f_{N}(B) \\
&=&\cup _{n=1}^{N}f_{n}(B),\text{ for each }B\in \mathcal{C}(X)
\end{eqnarray*}%
and%
\begin{eqnarray*}
S(B) &=&g_{1}(B)\cup g_{2}(B)\cup \cdot \cdot \cdot \cup g_{N}(B) \\
&=&\cup _{n=1}^{N}g_{n}(B),\text{ for each }B\in \mathcal{C}(X)
\end{eqnarray*}%
is generalized rational contractive. Then $T$ and $S$ have a unique common
attractor $A\in \mathcal{C}\left( X\right) ,$ that is,%
\begin{equation*}
A=T\left( A\right) =\cup _{n=1}^{N}f_{n}(A)=S\left( A\right) =\cup
_{n=1}^{N}g_{n}(A).
\end{equation*}%
Furthermore, for the initial set $B_{0}\in \mathcal{C}\left( X\right) $, the
iterative sequence of compact sets defined as%
\begin{equation*}
B_{2n+1}=T\left( B_{2n}\right) \text{ and }B_{2n+2}=S\left( B_{2n+1}\right) 
\text{ for }n=0,1,2,...
\end{equation*}%
converges to the common attractor of $T$ and $S$.

\noindent \textbf{Proof}\textit{.} \  \ Let $B_{0}\ $be any arbitrary element
in $\mathcal{C}\left( X\right) .$ Define%
\begin{equation*}
B_{1}=T(B_{0}),B_{2}=S(B_{1})
\end{equation*}%
and in general, we have a sequence $\left \{ B_{n}\right \} $ defined as $%
B_{2n+1}=T\left( B_{2n}\right) $ and $B_{2n+2}=S\left( B_{2n+1}\right) $ for 
$n=0,1,2,..._{.}$

\noindent Now, as the pair $\left( T,S\right) $ is \ a generalized rational
contraction, we have%
\begin{eqnarray*}
&&H_{\delta }(B_{2n+1},B_{2n+2})=H_{\delta }(T\left( B_{2n}\right) ,S\left(
B_{2n+1}\right) )\leq \alpha M_{T,S}(B_{2n},B_{2n+1}) \\
&=&\alpha \max \{H_{\delta }(B_{2n},B_{2n+1}),H_{\delta
}(B_{2n},T(B_{2n})),H_{\delta }(B_{2n+1},S(B_{2n+1})), \\
&&\frac{H_{\delta }(B_{2n},S\left( B_{2n+1}\right) )+H_{\delta
}(B_{2n+1},T\left( B_{2n}\right) )}{2}, \\
&&\frac{H_{\delta }\left( B_{2n+1},S\left( B_{2n+1}\right) \right)
[1+H_{\delta }\left( B_{2n+1},T\left( B_{2n}\right) \right) ]}{1+H_{\delta
}(B_{2n},B_{2n+1})}, \\
&&\frac{H_{\delta }(B_{2n+1},S\left( B_{2n+1}\right) [1+H_{\delta }\left(
B_{2n},T\left( B_{2n}\right) \right) ]}{1+H_{\delta }(B_{2n},B_{2n+1})}, \\
&&\frac{H_{\delta }(B_{2n+1},T\left( B_{2n}\right) )[1+H_{\delta
}(B_{2n+1},T\left( B_{2n}\right) )]}{1+H_{\delta }(B_{2n},B_{2n+1})}, \\
&&\frac{H_{\delta }(B_{2n+1},T\left( B_{2n}\right) )[1+H_{\delta
}(B_{2n},T\left( B_{2n}\right) )]}{1+H_{\delta }(B_{2n},B_{2n+1})}\}
\end{eqnarray*}%
\begin{eqnarray*}
&=&\alpha \max \{H_{\delta }(B_{2n},B_{2n+1}),H_{\delta
}(B_{2n},B_{2n+1}),H_{\delta }(B_{2n+1},B_{2n+2}), \\
&&\frac{H_{\delta }(B_{2n},B_{2n+2})+H_{\delta }(B_{2n+1},B_{2n+1})}{2}, \\
&&\frac{H_{\delta }\left( B_{2n+1},B_{2n+2}\right) [1+H_{\delta }\left(
B_{2n+1},B_{2n+1}\right) ]}{1+H_{\delta }(B_{2n},B_{2n+1})}, \\
&&\frac{H_{\delta }(B_{2n+1},B_{2n+2})[1+H_{\delta }\left(
B_{2n},B_{2n+1}\right) ]}{1+H_{\delta }(B_{2n},B_{2n+1})}, \\
&&\frac{H_{\delta }(B_{2n+1},B_{2n+1})[1+H_{\delta }(B_{2n+1},B_{2n+1})]}{%
1+H_{\delta }(B_{2n},B_{2n+1})}, \\
&&\frac{H_{\delta }(B_{2n+1},B_{2n+1})[1+H_{\delta }(B_{2n},B_{2n+1})]}{%
1+H_{\delta }(B_{2n},B_{2n+1})}\} \\
&=&\alpha \max \{H_{\delta }(B_{2n},B_{2n+1}),H_{\delta
}(B_{2n+1},B_{2n+2})\} \\
&=&\alpha H_{\delta }\left( B_{2n},B_{2n+1}\right) 
\end{eqnarray*}%
for all $n\in \{0,1,2,...\}.$

\noindent Similarly, we have%
\begin{eqnarray*}
&&H_{\delta }(B_{2n+2},B_{2n+3})=H_{\delta }(T\left( B_{2n+1}\right)
,S\left( B_{2n+2}\right) )=H_{\delta }(T\left( B_{2n+1}\right) ,S\left(
B_{2n+2}\right) ) \\
&\leq &\alpha M_{T,S}(B_{2n+1},B_{2n+2})
\end{eqnarray*}%
\begin{eqnarray*}
&=&\alpha \max \{H_{\delta }(B_{2n+1},B_{2n+2}),H_{\delta
}(B_{2n+1},T(B_{2n+1})),H_{\delta }(B_{2n+2},S(B_{2n+2})), \\
&&\frac{H_{\delta }(B_{2n+1},S\left( B_{2n+2}\right) )+H_{\delta
}(B_{2n+2},T\left( B_{2n+1}\right) )}{2}, \\
&&\frac{H_{\delta }\left( B_{2n+2},S\left( B_{2n+1}\right) \right)
[1+H_{\delta }\left( B_{2n+2},T\left( B_{2n+1}\right) \right) ]}{1+H_{\delta
}(B_{2n},B_{2n+1})}, \\
&&\frac{H_{\delta }(B_{2n+2},S\left( B_{2n+2}\right) [1+H_{\delta }\left(
B_{2n+1},T\left( B_{2n+1}\right) \right) ]}{1+H_{\delta }(B_{2n},B_{2n+1})},
\\
&&\frac{H_{\delta }(B_{2n+2},T\left( B_{2n+1}\right) )[1+H_{\delta
}(B_{2n+2},T\left( B_{2n+1}\right) )]}{1+H_{\delta }(B_{2n},B_{2n+1})}, \\
&&\frac{H_{\delta }(B_{2n+2},T\left( B_{2n+1}\right) )[1+H_{\delta
}(B_{2n+1},T\left( B_{2n+1}\right) )]}{1+H_{\delta }(B_{2n},B_{2n+1})}\}
\end{eqnarray*}%
\begin{eqnarray*}
&=&\alpha \max \{H_{\delta }(B_{2n+1},B_{2n+2}),H_{\delta
}(B_{2n+1},B_{2n+2}),H_{\delta }(B_{2n+2},B_{2n+3}), \\
&&\frac{H_{\delta }(B_{2n+1},B_{2n+3})+H_{\delta }(B_{2n+2},B_{2n+2})}{2}, \\
&&\frac{H_{\delta }\left( B_{2n+2},B_{2n+2}\right) [1+H_{\delta }\left(
B_{2n+2},B_{2n+2}\right) ]}{1+H_{\delta }(B_{2n+2},B_{2n+3})}, \\
&&\frac{H_{\delta }(B_{2n+2},B_{2n+3})[1+H_{\delta }\left(
B_{2n+1},B_{2n+2}\right) ]}{1+H_{\delta }(B_{2n+2},B_{2n+3})}, \\
&&\frac{H_{\delta }(B_{2n+1},B_{2n+2})[1+H_{\delta }(B_{2n+1},B_{2n+2})]}{%
1+H_{\delta }(B_{2n},B_{2n+1})}, \\
&&\frac{H_{\delta }(B_{2n+2},B_{2n+2})[1+H_{\delta }(B_{2n+1},B_{2n+2})]}{%
1+H_{\delta }(B_{2n+2},B_{2n+3})}\} \\
&=&\alpha \max \{H_{\delta }(B_{2n+1},B_{2n+2}),H_{\delta
}(B_{2n+2},B_{2n+3})\} \\
&=&\alpha H_{\delta }(B_{2n+1},B_{2n+2}).
\end{eqnarray*}%
Thus, we get%
\begin{equation*}
H_{\delta }(B_{2n+2},B_{2n+3})\leq \alpha H_{\delta }(B_{2n+1},B_{2n+2}).
\end{equation*}%
for all $n\in \{0,1,2,...\}.$ By continuation of this way, we obtain%
\begin{eqnarray}
H_{\delta }(B_{n},B_{n+1}) &\leq &\alpha H_{\delta }(B_{n-1},B_{n})  \notag
\\
&\leq &\alpha ^{2}H_{\delta }(B_{n-2},B_{n-1})  \notag \\
&\leq &...  \notag \\
&\leq &\alpha ^{n}H_{\delta }(B_{0},B_{1}). 
\end{eqnarray}%
Now for $m,n\in 
\mathbb{N}
$ with $m>n,$ we have%
\begin{eqnarray*}
H_{\delta }(B_{n},B_{m}) &\leq &H_{\delta }(B_{n},B_{n+1})+H_{\delta
}(B_{n+1},B_{n+2})+...+H_{\delta }(B_{m-1},B_{m}) \\
&\leq &\alpha ^{n}H_{\delta }(B_{0},B_{1})+\alpha ^{n+1}H_{\delta
}(B_{0},B_{1})+...+\alpha ^{m-1}H_{\delta }(B_{0},B_{1}) \\
&=&(\alpha ^{n}+\alpha ^{n+1}+...+\alpha ^{m-1})H_{\delta }(B_{0},B_{1}) \\
&\leq &\alpha ^{n}(1+\alpha +\alpha ^{2}+...)H_{\delta }(B_{0},B_{1}) \\
&=&\frac{\alpha ^{n}}{1-\alpha }H_{\delta }(B_{0},B_{1}).
\end{eqnarray*}%
It follows that by taking the limits as $n\rightarrow \infty ,$ we have $%
H_{\delta }(B_{n},B_{m})\rightarrow 0,$ that is, the sequence $\{B_{n}\}$ is
a Cauchy. Since $(C(X),H_{\delta })$ is complete dislocated metric space,
there exists a set $A^{\ast }\in C(X)$ such that $B_{n}$ $\rightarrow
A^{\ast }$as $n\rightarrow \infty ;$ that is, $\underset{n\rightarrow \infty 
}{\lim }H_{\delta }(B_{n},A^{\ast })=\underset{n\rightarrow \infty }{\lim }%
H_{\delta }(B_{n},B_{n+1})=H_{\delta }(A^{\ast },A^{\ast }).$

It follows from (1), that $H_{\delta }(B_{n},B_{n+1})\leq \alpha
^{n}H_{\delta }(B_{0},B_{1})\ $and by taking the limit as $n\rightarrow
\infty $ implies that $\underset{n\rightarrow \infty }{\lim }H_{\delta
}(B_{n},B_{n+1})=0$ and due to uniqueness of limit, we obtain $H_{\delta
}(A^{\ast },A^{\ast })=0.$

\noindent Now we are to show that $A^{\ast }$ is the common attractor of $T$
and $S$. Now%
\begin{eqnarray}
H_{\delta }(A^{\ast },S(A^{\ast })) &\leq &H_{\delta }(A^{\ast },T\left(
B_{2n}\right) )+H_{\delta }(T\left( B_{2n}\right) ,S(A^{\ast }))  \notag \\
&\leq &H_{\delta }(A^{\ast },B_{2n+1})+\alpha M_{T,S}\left( B_{2n},A^{\ast
}\right) ,  
\end{eqnarray}%
where%
\begin{eqnarray*}
M_{T,S}(B_{2n},A^{\ast }) &=&\alpha \max \{H_{\delta }(B_{2n},A^{\ast
}),H_{\delta }(B_{2n},T(B_{2n})),H_{\delta }(A^{\ast },S(A^{\ast })) \\
&&\frac{H_{\delta }(B_{2n},S\left( A^{\ast }\right) )+H_{\delta }(A^{\ast
},T\left( B_{2n}\right) )}{2}, \\
&&\frac{H_{\delta }\left( A^{\ast },S\left( B_{2n}\right) \right)
[1+H_{\delta }\left( A^{\ast },T\left( B_{2n}\right) \right) ]}{1+H_{\delta
}(B_{2n},A^{\ast })}, \\
&&\frac{H_{\delta }(A^{\ast },S\left( A^{\ast }\right) [1+H_{\delta }\left(
B_{2n},T\left( B_{2n}\right) \right) ]}{1+H_{\delta }(B_{2n},A^{\ast })}, \\
&&\frac{H_{\delta }(A^{\ast },T\left( B_{2n}\right) )[1+H_{\delta }(A^{\ast
},T\left( B_{2n}\right) )]}{1+H_{\delta }(B_{2n},A^{\ast })}, \\
&&\frac{H_{\delta }(A^{\ast },T\left( B_{2n}\right) )[1+H_{\delta
}(B_{2n},T\left( B_{2n}\right) )]}{1+H_{\delta }(B_{2n},A^{\ast })}\}. \\
&=&\alpha \max \{H_{\delta }(B_{2n},A^{\ast }),H_{\delta
}(B_{2n},B_{2n+1}),H_{\delta }(A^{\ast },S(A^{\ast })), \\
&&\frac{H_{\delta }(B_{2n},S\left( A^{\ast }\right) )+H_{\delta }(A^{\ast
},B_{2n+1})}{2}, \\
&&\frac{H_{\delta }\left( A^{\ast },B_{2n+1}\right) [1+H_{\delta }\left(
A^{\ast },B_{2n+1}\right) ]}{1+H_{\delta }(B_{2n},A^{\ast })}, \\
&&\frac{H_{\delta }(A^{\ast },S\left( A^{\ast }\right) [1+H_{\delta }\left(
B_{2n},B_{2n+1}\right) ]}{1+H_{\delta }(B_{2n},A^{\ast })}, \\
&&\frac{H_{\delta }(A^{\ast },B_{2n+1})[1+H_{\delta }(A^{\ast },B_{2n+1})]}{%
1+H_{\delta }(B_{2n},A^{\ast })}, \\
&&\frac{H_{\delta }(A^{\ast },B_{2n+1})[1+H_{\delta }(B_{2n},B_{2n+1})]}{%
1+H_{\delta }(B_{2n},A^{\ast })}\}.
\end{eqnarray*}%
Now, we have the following cases:

\begin{description}
\item[(i)] if $M_{T,S}\left( B_{2n},A^{\ast }\right) =H_{\delta }\left(
B_{2n},A^{\ast }\right) ,$ then by (3.2), we obtain%
\begin{equation*}
H_{\delta }(A^{\ast },S(A^{\ast }))\leq H_{\delta }(A^{\ast
},B_{2n+1})+\alpha H_{\delta }\left( B_{2n},A^{\ast }\right)
\end{equation*}%
and on taking limit as $n\rightarrow \infty $ implies%
\begin{eqnarray*}
\lim_{n\rightarrow \infty }H_{\delta }(A^{\ast },S(A^{\ast })) &\leq
&\lim_{n\rightarrow \infty }H_{\delta }(A^{\ast },B_{2n+1})+\alpha
\lim_{n\rightarrow \infty }H_{\delta }\left( B_{2n},A^{\ast }\right) \\
&=&0
\end{eqnarray*}%
and we obtain that $A^{\ast }=S\left( A^{\ast }\right) $.

\item[(ii)] In case $M_{T,S}\left( B_{2n},A^{\ast }\right) =H_{\delta
}\left( B_{2n},B_{2n+1}\right) ,$ then by (3.2), we have%
\begin{equation*}
H_{\delta }(A^{\ast },S(A^{\ast }))\leq H_{\delta }(A^{\ast
},B_{2n+1})+\alpha H\left( B_{2n},B_{2n+1}\right)
\end{equation*}%
and on taking limit as $n\rightarrow \infty $ implies%
\begin{eqnarray*}
\lim_{n\rightarrow \infty }H_{\delta }(A^{\ast },S(A^{\ast })) &\leq
&\lim_{n\rightarrow \infty }H_{\delta }(A^{\ast },B_{2n+1})+\alpha
\lim_{n\rightarrow \infty }H_{\delta }\left( B_{2n},B_{2n+1}\right) \\
&=&0
\end{eqnarray*}%
which implies $A^{\ast }=S\left( A^{\ast }\right) $.

\item[(iii)] In case $M_{T,S}\left( B_{2n},A^{\ast }\right) =H_{\delta
}\left( A^{\ast },S(A^{\ast })\right) ,$ then by (3.2), we obtain%
\begin{equation*}
H_{\delta }(A^{\ast },S(A^{\ast }))\leq H_{\delta }(A^{\ast
},B_{2n+1})+\alpha H_{\delta }\left( A^{\ast },S(A^{\ast })\right)
\end{equation*}%
and on taking limit as $n\rightarrow \infty $ implies%
\begin{equation*}
H_{\delta }(A^{\ast },S(A^{\ast }))\leq \alpha H_{\delta }\left( A^{\ast
},S(A^{\ast })\right) ,
\end{equation*}%
which implies%
\begin{equation*}
(1-\alpha )H_{\delta }\left( A^{\ast },S(A^{\ast })\right) \leq 0
\end{equation*}%
and since $1-\alpha >0,$ so we get $A^{\ast }=S\left( A^{\ast }\right) $.

\item[(iv)] In case $M_{T,S}\left( B_{2n},A^{\ast }\right) =\dfrac{H_{\delta
}(B_{2n},S\left( A^{\ast }\right) )+H_{\delta }(A^{\ast },B_{2n+1})}{2},$then%
\begin{eqnarray*}
H_{\delta }(A^{\ast },S(A^{\ast })) &\leq &H_{\delta }(A^{\ast },B_{2n+1})+%
\frac{\alpha }{2}[H_{\delta }(B_{2n},S\left( A^{\ast }\right) )+H_{\delta
}(A^{\ast },B_{2n+1})] \\
&\leq &H_{\delta }(A^{\ast },B_{2n+1})+\frac{\alpha }{2}[H_{\delta
}(B_{2n},A^{\ast })+H_{\delta }(A^{\ast },S\left( A^{\ast }\right)
)+H_{\delta }(A^{\ast },B_{2n+1})]
\end{eqnarray*}%
and on taking limit as $n\rightarrow \infty $ implies%
\begin{equation*}
H_{\delta }(A^{\ast },S(A^{\ast }))\leq \frac{\alpha }{2}H_{\delta }\left(
A^{\ast },S(A^{\ast })\right) ,
\end{equation*}%
which implies%
\begin{equation*}
(1-\frac{\alpha }{2})H_{\delta }\left( A^{\ast },S(A^{\ast })\right) \leq 0
\end{equation*}%
and since $1-\dfrac{\alpha }{2}>0,$ so we get $A^{\ast }=S\left( A^{\ast
}\right) $.

\item[(v)] In case $M_{T,S}\left( B_{2n},A^{\ast }\right) =\dfrac{H_{\delta
}\left( A^{\ast },B_{2n+1}\right) [1+H_{\delta }\left( A^{\ast
},B_{2n+1}\right) ]}{1+H_{\delta }(B_{2n},A^{\ast })},$ then%
\begin{equation*}
H_{\delta }(A^{\ast },S(A^{\ast }))\leq H_{\delta }(A^{\ast
},B_{2n+1})+\alpha \lbrack \frac{H_{\delta }\left( A^{\ast },B_{2n+1}\right)
[1+H_{\delta }\left( A^{\ast },B_{2n+1}\right) ]}{1+H_{\delta
}(B_{2n},A^{\ast })}].
\end{equation*}%
And on taking limit as $n\rightarrow \infty $ implies%
\begin{equation*}
H_{\delta }(A^{\ast },S(A^{\ast }))=0,
\end{equation*}%
which gives $A^{\ast }=S(A^{\ast })$.

\item[(vi)] When $M_{T,S}\left( B_{2n},A^{\ast }\right) =\dfrac{H_{\delta
}(A^{\ast },S\left( A^{\ast }\right) )[1+H_{\delta }(B_{2n},B_{2n+1})]}{%
1+H_{\delta }(B_{2n},A^{\ast })},$ we obtain%
\begin{equation*}
H_{\delta }(A^{\ast },S(A^{\ast }))\leq H_{\delta }(A^{\ast
},B_{2n+1})+\alpha \lbrack \frac{H_{\delta }(A^{\ast },S\left( A^{\ast
}\right) )[1+H_{\delta }(B_{2n},B_{2n+1})]}{1+H_{\delta }(B_{2n},A^{\ast })}%
].
\end{equation*}%
and on taking limit as $n\rightarrow \infty $ implies%
\begin{equation*}
H_{\delta }(A^{\ast },S(A^{\ast }))\leq \alpha H_{\delta }(A^{\ast },S\left(
A^{\ast }\right) ),
\end{equation*}%
since $1-\alpha >0,$ we get $A^{\ast }=S(A^{\ast })$.

\item[(vii)] In case $M_{T,S}\left( B_{2n},A^{\ast }\right) =\dfrac{%
H_{\delta }(A^{\ast },B_{2n+1})[1+H_{\delta }(A^{\ast },B_{2n+1})]}{%
1+H_{\delta }(B_{2n},A^{\ast })},$ then%
\begin{equation*}
H_{\delta }(A^{\ast },S(A^{\ast }))\leq H_{\delta }(A^{\ast
},B_{2n+1})+\alpha \lbrack \frac{H_{\delta }\left( A^{\ast },B_{2n+1}\right)
[1+H_{\delta }\left( A^{\ast },B_{2n+1}\right) ]}{1+H_{\delta
}(B_{2n},A^{\ast })}]
\end{equation*}%
and on taking limit as $n\rightarrow \infty $ implies%
\begin{equation*}
H_{\delta }(A^{\ast },S(A^{\ast }))=0
\end{equation*}%
which gives $A^{\ast }=S(A^{\ast })$.

\item[(viii)] Finally if $M_{T,S}\left( B_{2n},A^{\ast }\right) =\dfrac{%
H_{\delta }(A^{\ast },B_{2n+1})[1+H_{\delta }(B_{2n},B_{2n+1})]}{1+H_{\delta
}(B_{2n},A^{\ast })},$ then%
\begin{equation*}
H_{\delta }(A^{\ast },S(A^{\ast }))\leq H_{\delta }(A^{\ast
},B_{2n+1})+\alpha \lbrack \frac{H_{\delta }(A^{\ast },B_{2n+1})[1+H_{\delta
}(B_{2n},B_{2n+1})]}{1+H_{\delta }(B_{2n},A^{\ast })}]
\end{equation*}%
and on taking limit as $n\rightarrow \infty $ implies%
\begin{equation*}
H_{\delta }\left( A^{\ast },S\left( A^{\ast }\right) \right) =0
\end{equation*}%
and so $A^{\ast }=S\left( A^{\ast }\right) .$
\end{description}

\noindent Thus from all cases, we obtain that $A^{\ast }=S\left( A^{\ast
}\right) .$

\noindent Again, we have%
\begin{eqnarray}
H_{\delta }(A^{\ast },T(A^{\ast })) &\leq &H_{\delta }(A^{\ast },S\left(
B_{2n+1}\right) )+H_{\delta }(T(A^{\ast }),S\left( B_{2n+1}\right) )  \notag
\\
&\leq &H_{\delta }(A^{\ast },B_{2n+2})+\alpha M_{T,S}\left( A^{\ast
},B_{2n+1}\right) ,  
\end{eqnarray}%
where%
\begin{eqnarray*}
M_{T,S}(A^{\ast },B_{2n+1}) &=&\alpha \max \{H_{\delta }(A^{\ast
},B_{2n+1}),H_{\delta }(A^{\ast },T(A^{\ast })),H_{\delta
}(B_{2n+1},S(B_{2n+1})), \\
&&\frac{H_{\delta }(A^{\ast },S\left( B_{2n+1}\right) )+H_{\delta
}(B_{2n+1},T\left( A^{\ast }\right) )}{2}, \\
&&\frac{H_{\delta }\left( B_{2n+1},S\left( A^{\ast }\right) \right)
[1+H_{\delta }\left( B_{2n+1},T\left( A^{\ast }\right) \right) ]}{%
1+H_{\delta }(A^{\ast },B_{2n+1})}, \\
&&\frac{H_{\delta }(B_{2n+1},S\left( B_{2n+1}\right) [1+H_{\delta }\left(
A^{\ast },T\left( A^{\ast }\right) \right) ]}{1+H_{\delta }(A^{\ast
},B_{2n+1})}, \\
&&\frac{H_{\delta }(B_{2n+1},T\left( A^{\ast }\right) )[1+H_{\delta
}(B_{2n+1},T\left( A^{\ast }\right) )]}{1+H_{\delta }(A^{\ast },B_{2n+1})},
\\
&&\frac{H_{\delta }(B_{2n+1},T\left( A^{\ast }\right) )[1+H_{\delta
}(A^{\ast },T\left( A^{\ast }\right) )]}{1+H_{\delta }(A^{\ast },B_{2n+1})}\}
\\
&=&\alpha \max \{H_{\delta }(A^{\ast },B_{2n+1}),H_{\delta }(A^{\ast
},T(A^{\ast })),H_{\delta }(B_{2n+1},B_{2n+2}), \\
&&\frac{[H_{\delta }(A^{\ast },B_{2n+2})+H_{\delta }(B_{2n+1},T\left(
A^{\ast }\right) )]}{2}, \\
&&\frac{H_{\delta }\left( B_{2n+1},A^{\ast }\right) [1+H_{\delta }\left(
B_{2n+1},T\left( A^{\ast }\right) \right) ]}{1+H_{\delta }(A^{\ast
},B_{2n+1})}, \\
&&\frac{H_{\delta }(B_{2n+1},B_{2n+2})[1+H_{\delta }\left( A^{\ast },T\left(
A^{\ast }\right) \right) ]}{1+H_{\delta }(A^{\ast },B_{2n+1})}, \\
&&\frac{H_{\delta }(B_{2n+1},T\left( A^{\ast }\right) )[1+H_{\delta
}(B_{2n+1},T\left( A^{\ast }\right) )]}{1+H_{\delta }(A^{\ast },B_{2n+1})},
\\
&&\frac{H_{\delta }(B_{2n+1},T\left( A^{\ast }\right) )[1+H_{\delta
}(A^{\ast },T\left( A^{\ast }\right) )]}{1+H_{\delta }(A^{\ast },B_{2n+1})}%
\}.
\end{eqnarray*}%
Now, we have again the following cases:

\begin{description}
\item[(i)] If $M_{T,S}\left( B_{2n+1},A^{\ast }\right) =H_{\delta }(A^{\ast
},B_{2n+1}),$ then by (3), we have%
\begin{equation*}
H_{\delta }(A^{\ast },T(A^{\ast }))\leq H_{\delta }(A^{\ast
},B_{2n+2})+\alpha H_{\delta }(A^{\ast },B_{2n+1})
\end{equation*}%
and on taking limit as $n\rightarrow \infty $ implies%
\begin{equation*}
H_{\delta }(A^{\ast },T(A^{\ast }))=0,
\end{equation*}%
which implies $A^{\ast }=T\left( A^{\ast }\right) $.

\item[(ii)] If $M_{T,S}\left( B_{2n+1},A^{\ast }\right) =H_{\delta }(A^{\ast
},T(A^{\ast })),$ then by (3), we get%
\begin{equation*}
H_{\delta }(A^{\ast },T(A^{\ast }))\leq H_{\delta }(A^{\ast
},B_{2n+2})+\alpha H_{\delta }(A^{\ast },T(A^{\ast }))
\end{equation*}%
and on taking limit as $n\rightarrow \infty $ implies%
\begin{equation*}
H_{\delta }(A^{\ast },T(A^{\ast }))\leq 0+\alpha H_{\delta }(A^{\ast
},T(A^{\ast })),
\end{equation*}%
that is,%
\begin{equation*}
(1-\alpha )H_{\delta }(A^{\ast },T(A^{\ast }))\leq 0
\end{equation*}%
and since $(1-\alpha )>0$, so we get $A^{\ast }=T\left( A^{\ast }\right) .$

\item[(iii)] If $M_{T,S}\left( B_{2n+1},A^{\ast }\right) =H_{\delta }\left(
B_{2n+1},B_{2n+2}\right) ,$ then by (3), we obtain that%
\begin{equation*}
H_{\delta }(A^{\ast },T(A^{\ast }))\leq H_{\delta }(A^{\ast
},B_{2n+2})+\alpha H_{\delta }\left( B_{2n+1},B_{2n+2}\right) 
\end{equation*}%
and on taking limit as $n\rightarrow \infty $ implies%
\begin{equation*}
H_{\delta }(A^{\ast },T(A^{\ast }))=0,
\end{equation*}%
that is, $A^{\ast }=T\left( A^{\ast }\right) $.

\item[(iv)] If $M_{T,S}\left( B_{2n+1},A^{\ast }\right) =\dfrac{H_{\delta
}(A^{\ast },B_{2n+2})+H_{\delta }(B_{2n+1},T\left( A^{\ast }\right) )}{2},$
then%
\begin{equation*}
H_{\delta }(A^{\ast },T(A^{\ast }))\leq H_{\delta }(A^{\ast },B_{2n+2})+%
\frac{\alpha }{2}[H_{\delta }(A^{\ast },B_{2n+2})+H_{\delta
}(B_{2n+1},T\left( A^{\ast }\right) )]
\end{equation*}%
and on taking limit as $n\rightarrow \infty $ implies%
\begin{equation*}
H_{\delta }(A^{\ast },T(A^{\ast }))\leq 0+\frac{\alpha }{2}H_{\delta
}(A^{\ast },T(A^{\ast })),
\end{equation*}%
that is,%
\begin{equation*}
(1-\frac{\alpha }{2})H_{\delta }(A^{\ast },T(A^{\ast }))\leq 0
\end{equation*}%
and since $1-\dfrac{\alpha }{2}>0$, so we get $A^{\ast }=T\left( A^{\ast
}\right) .$

\item[(v)] If $M_{T,S}\left( B_{2n+1},A^{\ast }\right) =\dfrac{H_{\delta
}\left( B_{2n+1},A^{\ast }\right) [1+H_{\delta }\left( B_{2n+1},T\left(
A^{\ast }\right) \right) ]}{1+H_{\delta }(A^{\ast },B_{2n+1})},$ we have%
\begin{equation*}
H_{\delta }(A^{\ast },T(A^{\ast }))\leq H_{\delta }(A^{\ast
},B_{2n+2})+\alpha \frac{H_{\delta }\left( B_{2n+1},A^{\ast }\right)
[1+H_{\delta }\left( B_{2n+1},T\left( A^{\ast }\right) \right) ]}{%
1+H_{\delta }(A^{\ast },B_{2n+1})}
\end{equation*}%
and on taking limit as $n\rightarrow \infty $ implies%
\begin{equation*}
H_{\delta }(A^{\ast },T(A^{\ast }))=0,
\end{equation*}%
that is, $A^{\ast }=T\left( A^{\ast }\right) $.

\item[(vi)] If $M_{T,S}\left( B_{2n+1},A^{\ast }\right) =\dfrac{H_{\delta
}(B_{2n+1},B_{2n+2})(1+H_{\delta }\left( A^{\ast },T\left( A^{\ast }\right)
\right) }{1+H_{\delta }(A^{\ast },B_{2n+1})},$ then%
\begin{equation*}
H_{\delta }(A^{\ast },T(A^{\ast }))\leq H_{\delta }(A^{\ast
},B_{2n+2})+\alpha \frac{H_{\delta }(B_{2n+1},B_{2n+2})(1+H_{\delta }\left(
A^{\ast },T\left( A^{\ast }\right) \right) }{1+H_{\delta }(A^{\ast
},B_{2n+1})}
\end{equation*}%
and on taking limit as $n\rightarrow \infty $ implies%
\begin{equation*}
H_{\delta }(A^{\ast },T(A^{\ast }))=0,
\end{equation*}%
that is, $A^{\ast }=T\left( A^{\ast }\right) .$

\item[(vii)] If $M_{T,S}\left( B_{2n+1},A^{\ast }\right) =\dfrac{H_{\delta
}(B_{2n+1},T\left( A^{\ast }\right) )[1+H_{\delta }(B_{2n+1},T\left( A^{\ast
}\right) )]}{1+H_{\delta }(A^{\ast },B_{2n+1})},$%
\begin{eqnarray*}
&&H_{\delta }(A^{\ast },T(A^{\ast }))\leq H_{\delta }(A^{\ast
},B_{2n+2})+\alpha \frac{H_{\delta }(B_{2n+1},T\left( A^{\ast }\right)
)[1+H_{\delta }(B_{2n+1},T\left( A^{\ast }\right) )]}{1+H_{\delta }(A^{\ast
},B_{2n+1})} \\
&\leq &H_{\delta }(A^{\ast },B_{2n+2})+\alpha \frac{\lbrack H_{\delta
}(B_{2n+1},A^{\ast })+H_{\delta }(A^{\ast },T\left( A^{\ast }\right)
)][1+H_{\delta }(B_{2n+1},T\left( A^{\ast }\right) )]}{1+H_{\delta }(A^{\ast
},B_{2n+1})}
\end{eqnarray*}%
and on taking limit as $n\rightarrow \infty $ implies%
\begin{equation*}
H_{\delta }(A^{\ast },T(A^{\ast }))\leq \alpha H_{\delta }(A^{\ast
},T(A^{\ast })),
\end{equation*}%
that is,%
\begin{equation*}
(1-\alpha )H_{\delta }(A^{\ast },T(A^{\ast }))\leq 0
\end{equation*}%
and since $1-\alpha >0$, so we get $A^{\ast }=T\left( A^{\ast }\right) .$

\item[(viii)] Finally if $M_{T,S}\left( B_{2n},A^{\ast }\right) =\dfrac{%
H_{\delta }(B_{2n+1},T\left( A^{\ast }\right) )[1+H_{\delta }(A^{\ast
},T\left( A^{\ast }\right) )]}{1+H_{\delta }(A^{\ast },B_{2n+1})},$ then%
\begin{eqnarray*}
&&H_{\delta }(A^{\ast },T(A^{\ast }))\leq H_{\delta }(A^{\ast
},B_{2n+2})+\alpha \frac{H_{\delta }(B_{2n+1},T\left( A^{\ast }\right)
)[1+H_{\delta }(A^{\ast },T\left( A^{\ast }\right) )]}{1+H_{\delta }(A^{\ast
},B_{2n+1})} \\
&\leq &H_{\delta }(A^{\ast },B_{2n+2})+\alpha \frac{\lbrack H_{\delta
}(B_{2n+1},A^{\ast })+H_{\delta }(A^{\ast },T\left( A^{\ast }\right)
)][1+H_{\delta }(A^{\ast },T\left( A^{\ast }\right) )]}{1+H_{\delta
}(A^{\ast },B_{2n+1})}
\end{eqnarray*}%
and on taking limit as $n\rightarrow \infty $ implies%
\begin{equation*}
H_{\delta }\left( A^{\ast },T\left( A^{\ast }\right) \right) \leq \alpha
H_{\delta }(A^{\ast },T\left( A^{\ast }\right) ),
\end{equation*}%
and since $\left( 1-\alpha \right) >0$ which gives $A^{\ast }=T\left(
A^{\ast }\right) .$
\end{description}

\noindent Thus, $A^{\ast }$ is the common attractor of $T$ and $S.$

\noindent Finally, to prove that the common attractor of $T$ and $S$ is
unique, we suppose that $A^{\ast }$ and $B^{\ast }$ are the two common
attractors of $T$ and $S$ in $\mathcal{C}(X).$ Since the pair $(T,S)$ is
generalized rational contractive mappings, so we obtain that%
\begin{eqnarray*}
H_{\delta }(A^{\ast },B^{\ast }) &=&H_{\delta }(T(A^{\ast }),S(B^{\ast })) \\
&\leq &\alpha M_{T,S}\left( A^{\ast },B^{\ast }\right) ,
\end{eqnarray*}%
where%
\begin{eqnarray*}
M_{T,S}\left( A^{\ast },B^{\ast }\right) &=&\max \{H_{\delta }\left( A^{\ast
},B^{\ast }\right) ,H_{\delta }(A^{\ast },T\left( A^{\ast }\right)
),H_{\delta }(B^{\ast },S\left( B^{\ast }\right) ), \\
&&\frac{H_{\delta }(A^{\ast },S\left( B^{\ast }\right) )+H_{\delta }(B^{\ast
},T\left( A^{\ast }\right) )}{2},\dfrac{H_{\delta }(B^{\ast },S\left(
B^{\ast }\right) )[1+H_{\delta }(B^{\ast },T(A^{\ast }))]}{1+H_{\delta
}\left( A^{\ast },B^{\ast }\right) }, \\
&&\frac{H_{\delta }(B^{\ast },S\left( B^{\ast }\right) )[1+H_{\delta
}(A^{\ast },T\left( A^{\ast }\right) )]}{1+H_{\delta }\left( A^{\ast
},B^{\ast }\right) },\frac{H_{\delta }(B^{\ast },T\left( A^{\ast }\right)
)[1+H_{\delta }(B^{\ast },T\left( A^{\ast }\right) )]}{1+H_{\delta }\left(
A^{\ast },B^{\ast }\right) }, \\
&&\frac{H_{\delta }(B^{\ast },T\left( A^{\ast }\right) )[1+H_{\delta
}(A^{\ast },T\left( A^{\ast }\right) )]}{1+H_{\delta }\left( A^{\ast
},B^{\ast }\right) }\} \\
&=&\max \{H_{\delta }\left( A^{\ast },B^{\ast }\right) ,\frac{H_{\delta
}(A^{\ast },B^{\ast })+H_{\delta }(B^{\ast },A^{\ast })}{2},\frac{H_{\delta
}(B^{\ast },A^{\ast })[1+H_{\delta }(B^{\ast },A^{\ast })]}{1+H_{\delta
}\left( A^{\ast },B^{\ast }\right) }\}, \\
&=&H_{\delta }\left( A^{\ast },B^{\ast }\right) ,
\end{eqnarray*}%
that is,%
\begin{equation*}
H_{\delta }(A^{\ast },B^{\ast })\leq \alpha H_{\delta }\left( A^{\ast
},B^{\ast }\right) ,
\end{equation*}%
which implies $A^{\ast }=B^{\ast }$. Hence, $T$ and $S$ have a unique common
attractor $A^{\ast }$ in $\mathcal{C}(X)$. $\square $\medskip

\noindent \textbf{Remark 3.2.} \  \ In Theorem 3.1, let $\mathcal{S}(X)$ be
the class of all singleton subset of $X.$ Obviously $\mathcal{S}(X)\subseteq
C(X).$ Further, suppose that $f_{m}=f$ and $g_{m}=g$ for every $m,$ with $%
f=f_{1}$ and $g=g_{1}$, then the operator $T,S:\mathcal{S}(X)\rightarrow 
\mathcal{S}(X)$ will be%
\begin{equation*}
T(A)=f(A)\text{ for all }A\in \mathcal{S}(X)
\end{equation*}%
and%
\begin{equation*}
S(A)=g(A)\text{ for all }A\in \mathcal{S}(X).
\end{equation*}%
Follows from above setup, we obtain the following result.\medskip

\noindent \textbf{Corollary 3.3.} \  \ Let $(\mathcal{C}(X),\delta )$ be a
complete dislocated metric space and $f,g:X\rightarrow X$ be two mappings.\
Let $T,S:\mathcal{S}(X)\rightarrow \mathcal{S}(X)$ is the map defined above
in Remark 3.2 satisfying $\alpha \in \lbrack 0,1),$such that%
\begin{equation*}
H_{\delta }(T(U),S(V))\leq \alpha M_{T,S}(U,V),
\end{equation*}%
is satisfied for all $U,V\in C(X),$ where%
\begin{eqnarray*}
M_{T,S}(U,V) &=&\max \{H_{\delta }(U,V),H_{\delta }(U,T\left( U\right)
),H_{\delta }(V,S\left( V\right) ), \\
&&\frac{H_{\delta }(U,S\left( V\right) )+H_{\delta }(V,T\left( U\right) )}{2}%
,\dfrac{H_{\delta }(V,S\left( V\right) )[1+H_{\delta }(V,T(U))]}{1+H_{\delta
}\left( U,V\right) }, \\
&&\frac{H_{\delta }(V,S\left( V\right) )[1+H_{\delta }(U,T\left( U\right) )]%
}{1+H_{\delta }(U,V)},\frac{H_{\delta }(V,T\left( U\right) )[1+H_{\delta
}(V,T\left( U\right) )]}{1+H_{\delta }(U,V)}, \\
&&\frac{H_{\delta }(V,T\left( U\right) )[1+H_{\delta }(U,T\left( U\right) )]%
}{1+H_{\delta }(U,V)}\}.
\end{eqnarray*}%
Then $T$ and $S$ have at most one attractor, that is, there exists a unique $%
U^{\ast }\in \mathcal{S}(X)$ such that%
\begin{equation*}
U^{\ast }=T(U^{\ast })=S(U^{\ast }).
\end{equation*}%
Furthermore, for the singleton set $B_{0}\in \mathcal{S}(X),$ the iterative
sequence of compact sets defined as%
\begin{equation*}
B_{2n+1}=T\left( B_{2n}\right) ,\text{ }B_{2n+2}=S\left( B_{2n+1}\right) 
\text{ for }n=0,1,2,...
\end{equation*}%
converges to the common attractor of $T$ and $S$.\medskip

\noindent \textbf{Corollary 3.4. \  \ }Let $(\mathcal{C}(X),\delta )$ be a
complete dislocated metric space and $(X:f_{m},g_{m};m=1,2,$%
\textperiodcentered \textperiodcentered \textperiodcentered $,N)$ be
generalized iterated function system, where each $\left( f_{i},g_{i}\right)
\ $is a generalized contraction on $X$ for $i=1,2,...,N.\ $Then $T,S:%
\mathcal{C}(X)\rightarrow \mathcal{C}(X)$ defined in Theorem 3.1 has a
unique attractor in $\mathcal{C}\left( X\right) .$ Furthermore, for any set $%
B_{0}\in \mathcal{C}\left( X\right) $, the sequence of compact sets%
\begin{equation*}
B_{2n+1}=T\left( B_{2n}\right) ,\text{ }B_{2n+2}=S\left( B_{2n+1}\right) 
\text{ for }n=0,1,2,...
\end{equation*}%
converges to the common attractor of $T$ and $S$.\medskip

\noindent \textbf{Theorem 3.5} (\textbf{Generalized Collage Theorem)} Let $%
(X,\delta )$ be a dislocated metric space and $\left \{
X;f_{1},f_{2,}...,f_{N};g_{1},g_{2,}...,g_{N}\right \} $ be a given GIFS.
Suppose that the pair of self mappings $(T,S)$ defined by%
\begin{equation*}
T(B)=f_{1}(B)\cup f_{2}(B)\cup \cdot \cdot \cdot \cup f_{N}(B)\text{ for
each }B\in \mathcal{C}(X)
\end{equation*}%
and%
\begin{equation*}
S(B)=g_{1}(B)\cup g_{2}(B)\cup \cdot \cdot \cdot \cup g_{N}(B)\text{ for
each }B\in \mathcal{C}(X)
\end{equation*}%
is satisfying%
\begin{equation*}
H_{\delta }(T\left( A\right) ,S(A))\leq \alpha H_{\delta }(A,B)\text{ for }%
A,B\in \mathcal{C}(X),
\end{equation*}%
where $0\leq \alpha <1.$ If for any $A\in \mathcal{C}(X)$ and for $%
\varepsilon \geq 0$ such that either%
\begin{equation*}
H_{\delta }(A,T(A))\leq \varepsilon 
\end{equation*}%
or%
\begin{equation*}
H_{\delta }(A,S(A))\leq \varepsilon ,
\end{equation*}%
then%
\begin{equation*}
H_{\delta }(A,U)\leq \frac{\varepsilon }{1-\alpha },
\end{equation*}%
where $U\in C(X)$ is the common attractor of $T\,$and $S$.

\noindent \textbf{Proof.} \  \ It follows from Corollary 3.1 that $U$ is the
common attractor of mappings $T,S:C(X)\rightarrow C(X),$ that is, $U=T\left(
U\right) =S\left( U\right) .$

\noindent Also, for any $B_{0}$ in $\mathcal{C}\left( X\right) ,$ we define
a sequence $\left \{ B_{n}\right \} $ as $B_{2n+1}=T\left( B_{2n}\right) $
and $B_{2n+2}=S\left( B_{2n+1}\right) $ for $n=0,1,2,...,$ we have%
\begin{equation*}
\lim_{n\rightarrow \infty }H_{\delta }\left( T\left( B_{2n}\right) ,U\right)
=\lim_{n\rightarrow \infty }H_{\delta }\left( S\left( B_{2n+1}\right)
,U\right) =0.
\end{equation*}%
\noindent Now, as the pair $\left( T,S\right) $ is a generalized
contraction, we have%
\begin{equation*}
H_{\delta }\left( B_{n},B_{n+1}\right) \leq \alpha ^{n}H_{\delta }\left(
B_{0},B_{1}\right) .
\end{equation*}%
Assume that $H_{\delta }(A,T(A))\leq \varepsilon \ $for any $A\in \mathcal{C}%
(X).$\ Now, we have%
\begin{eqnarray*}
H_{\delta }(A,U) &\leq &H_{\delta }(A,T(A))+H_{\delta }(T(A),S(U)) \\
&\leq &\varepsilon +\alpha H_{\delta }(A,U),
\end{eqnarray*}%
which further implies that%
\begin{equation*}
H_{\delta }(A,U)\leq \frac{\varepsilon }{1-\alpha }.
\end{equation*}%
Similarly, if we assume that $H_{\delta }(A,S(A))\leq \varepsilon \ $for any 
$A\in \mathcal{C}(X).$\ Then, we have%
\begin{eqnarray*}
H_{\delta }(A,U) &\leq &H_{\delta }(A,S(A))+H_{\delta }(S(A),T(U)) \\
&\leq &\varepsilon +\alpha H_{\delta }(A,U),
\end{eqnarray*}%
which gives%
\begin{equation*}
H_{\delta }(A,U)\leq \frac{\varepsilon }{1-\alpha }.
\end{equation*}

\section{Well-posedness of Common Attractors}

In this section, we will define the well-posedness of attractors based
problems of rational contraction maps in the framework of dislocated metric
spaces.

\noindent \textbf{Definition 4.1. \  \ }Let $(X,\delta )$ be a dislocated
metric space. An attractor based problem of mapping $T:\mathcal{C}%
(X)\rightarrow \mathcal{C}(X)\ $is called well-posed if $T$ has a unique
attractor $A^{\ast }\in \mathcal{C}(X)$ and for any sequence $\{A_{n}\}$ in $%
\mathcal{C}(X)$ such that $\lim \limits_{n\rightarrow \infty }H_{\delta
}(T(A_{n}),A_{n})=0$ implies that $\lim \limits_{n\rightarrow \infty
}A_{n}=A^{\ast }.$\medskip

\noindent \textbf{Definition 4.2. \  \ }Let $(X,\delta )$ be a dislocated
metric space. A common attractor based problem of mappings $T,S:\mathcal{C}%
(X)\rightarrow \mathcal{C}(X)\ $is called well-posed if $T$ and $S$ have a
unique common attractor $A^{\ast }\in \mathcal{C}(X)$ and for any sequence $%
\{A_{n}\}$ in $\mathcal{C}(X)$ such that $\lim \limits_{n\rightarrow \infty
}H_{\delta }(T(A_{n}),A_{n})=0$ and\ $\lim \limits_{n\rightarrow \infty
}H_{\delta }(S(A_{n}),A_{n})=0$ implies that $\lim \limits_{n\rightarrow
\infty }A_{n}=A^{\ast }.$\medskip

\noindent \textbf{Theorem 4.3. \ } Let $(X,\delta )$ be a complete
dislocated metric space and $T$ and $S$ are self-maps on $\mathcal{C}(X)$ as
in Theorem 3.1. Then the common attractor based problem of $T$ and $S$ is
well-posed.

\noindent \textbf{Proof}. It follow from Theorem 3.1, that maps $T$ and $S$
have a unique common attractor say $B^{\ast }$ and $H_{\delta }\left(
B^{\ast },B^{\ast }\right) =0.$ Let $\{B_{n}\}$ be the sequence in $\mathcal{%
C}(X)$ such that $\lim \limits_{n\rightarrow \infty }H_{\delta
}(T(B_{n}),B_{n})=0$ and $\lim \limits_{n\rightarrow \infty }H_{\delta
}(S(B_{n}),B_{n})=0.$ We want to show that $B^{\ast }=\lim
\limits_{n\rightarrow \infty }B_{n}$ for every integer $n.$ As the pair of $%
\left( S,T\right) $ is generalized rational contractive operators, so that%
\begin{eqnarray}
H_{\delta }(B_{n},B^{\ast }) &=&H_{\delta }(B^{\ast },B_{n})  \notag \\
&\leq &H_{\delta }(B^{\ast },T(B_{2n}))+H_{\delta }(T(B_{2n}),B_{n})  \notag
\\
&=&H_{\delta }(S\left( B^{\ast }\right) ,T(B_{2n}))+H_{\delta
}(T(B_{2n}),B_{n})  \notag \\
&\leq &\alpha M_{T,S}(B^{\ast },B_{n})+H_{\delta }(T(B_{2n}),B_{n}). 
\end{eqnarray}%
where%
\begin{eqnarray*}
M_{T,S}(B^{\ast },B_{n}) &=&\max \{H_{\delta }(B^{\ast },B_{n}),H_{\delta
}(B^{\ast },T\left( B^{\ast }\right) ),H_{\delta }(B_{n},S\left(
B_{n}\right) ), \\
&&\frac{H_{\delta }(B^{\ast },S\left( B_{n}\right) )+H_{\delta
}(B_{n},T\left( B^{\ast }\right) )}{2}, \\
&&\dfrac{H_{\delta }(B_{n},S\left( B_{n}\right) )[1+H_{\delta
}(B_{n},T(B^{\ast }))]}{1+H_{\delta }(B^{\ast },B_{n})}, \\
&&\frac{H_{\delta }(B_{n},S\left( B_{n}\right) )[1+H_{\delta }(B^{\ast
},T\left( B^{\ast }\right) )]}{1+H_{\delta }(B^{\ast },B_{n})}, \\
&&\frac{H_{\delta }(B_{n},T\left( B^{\ast }\right) )[1+H_{\delta
}(B_{n},T\left( B^{\ast }\right) )]}{1+H_{\delta }(B^{\ast },B_{n})}, \\
&&\frac{H_{\delta }(B_{n},T\left( B^{\ast }\right) )[1+H_{\delta }(B^{\ast
},T\left( B^{\ast }\right) )]}{1+H_{\delta }(B^{\ast },B_{n})}\}.
\end{eqnarray*}%
Then following cases arise:

\begin{description}
\item[(i)] If $M_{T,S}(B^{\ast },B_{n})=H_{\delta }(B^{\ast },B_{n}),$ then%
\begin{equation*}
H_{\delta }(B_{n},B^{\ast })\leq \alpha H_{\delta }(B^{\ast
},B_{n})+H_{\delta }(T(B_{2n}),B_{n}),
\end{equation*}%
which implies%
\begin{equation*}
H_{\delta }(B^{\ast },B_{n})\leq \frac{1}{1-\alpha }H_{\delta
}(T(B_{2n}),B_{n})
\end{equation*}%
and on taking limit as $n\rightarrow \infty $ implies%
\begin{equation*}
\lim \limits_{n\rightarrow \infty }H_{\delta }(B^{\ast },B_{n})\leq 0,
\end{equation*}%
that is, $\lim \limits_{n\rightarrow \infty }B_{n}=B^{\ast }.$

\item[(ii)] If $M_{T,S}(B^{\ast },B_{n})=H_{\delta }(B^{\ast },T(B^{\ast })),
$ then%
\begin{eqnarray*}
H_{\delta }(B_{n},B^{\ast }) &\leq &\alpha H_{\delta }(B^{\ast },T(B^{\ast
}))+H_{\delta }(T(B_{2n}),B_{n}) \\
&\leq &H_{\delta }(T(B_{2n}),B_{n})
\end{eqnarray*}%
and on taking limit as $n\rightarrow \infty $ implies%
\begin{equation*}
\lim \limits_{n\rightarrow \infty }H(B^{\ast },B_{n})=0,
\end{equation*}%
which implies that $\lim \limits_{n\rightarrow \infty }B_{n}=B^{\ast }.$

\item[(iii)] If $M_{T,S}(B^{\ast },B_{n})=H_{\delta }(B_{n},S\left(
B_{n}\right) ),$ then%
\begin{equation*}
H_{\delta }(B_{n},B^{\ast })\leq \alpha H_{\delta }(B_{n},S\left(
B_{n}\right) )+H_{\delta }(T(B_{2n}),B_{n}),
\end{equation*}%
and on taking limit as $n\rightarrow \infty $ implies%
\begin{eqnarray*}
\lim \limits_{n\rightarrow \infty }H_{\delta }(B_{n},B^{\ast }) &\leq
&\lim \limits_{n\rightarrow \infty }[\alpha H_{\delta }(B_{n},S\left(
B_{n}\right) )+H_{\delta }(T(B_{2n}),B_{n}))] \\
&=&0,
\end{eqnarray*}%
which implies that $\lim \limits_{n\rightarrow \infty }B_{n}=B^{\ast }.$

\item[(iv)] If $M_{T,S}(B^{\ast },B_{n})=\dfrac{H_{\delta }(B^{\ast
},S\left( B_{n}\right) )+H_{\delta }(B_{n},T\left( B^{\ast }\right) )}{2}$
then%
\begin{eqnarray*}
H_{\delta }(B_{n},B^{\ast }) &\leq &\frac{\alpha }{2}[H_{\delta }(B^{\ast
},S\left( B_{n}\right) )+H_{\delta }(B_{n},B^{\ast })]+H_{\delta
}(T(B_{2n}),B_{n}) \\
&\leq &\frac{\alpha }{2}[H_{\delta }(B_{n},S\left( B_{n}\right) )+2H_{\delta
}(B_{n},B^{\ast })]+H_{\delta }(T(B_{2n}),B_{n}),
\end{eqnarray*}%
that is,%
\begin{equation*}
H_{\delta }(B_{n},B^{\ast })\leq \frac{\alpha }{2\left( 1-\alpha \right) }%
H_{\delta }(B_{n},S\left( B_{n}\right) )+\frac{1}{1-\alpha }H_{\delta
}(T(B_{2n}),B_{n})
\end{equation*}%
and on taking limit as $n\rightarrow \infty $ implies%
\begin{equation*}
\lim \limits_{n\rightarrow \infty }H_{\delta }(B_{n},B^{\ast })=0,
\end{equation*}%
which implies that $\lim \limits_{n\rightarrow \infty }B_{n}=B^{\ast }.$

\item[(v)] If $M_{T,S}(B^{\ast },B_{n})=\dfrac{H_{\delta }(B_{n},T\left(
B^{\ast }\right) )[1+H_{\delta }(B_{n},T\left( B^{\ast }\right) )]}{%
1+H_{\delta }(B^{\ast },B_{n})},$ then%
\begin{equation*}
H_{\delta }(B_{n},B^{\ast })\leq \alpha \left( \frac{H_{\delta
}(B_{n},B^{\ast })[1+H_{\delta }(B_{n},B^{\ast })]}{1+H_{\delta }(B^{\ast
},B_{n})}\right) +H_{\delta }(T(B_{2n}),B_{n}),
\end{equation*}%
On taking limit $n\rightarrow \infty $ implies%
\begin{eqnarray*}
\lim \limits_{n\rightarrow \infty }H_{\delta }(B_{n},B^{\ast }) &\leq
&\lim \limits_{n\rightarrow \infty }\alpha \left( \frac{H_{\delta
}(B_{n},B^{\ast })[1+H_{\delta }(B_{n},B^{\ast })]}{1+H_{\delta }(B^{\ast
},B_{n})}\right) +\lim \limits_{n\rightarrow \infty }H_{\delta
}(T(B_{2n}),B_{n}), \\
&=&0,
\end{eqnarray*}%
which implies that $\lim \limits_{n\rightarrow \infty }B_{n}=B^{\ast }.$

\item[(vi)] If $M_{T,S}(B^{\ast },B_{n})=\dfrac{H_{\delta }(B_{n},S\left(
B_{n}\right) )[1+H_{\delta }(B^{\ast },T\left( B^{\ast }\right) )]}{%
1+H_{\delta }(B^{\ast },B_{n})}$ then%
\begin{equation*}
H_{\delta }(B_{n},B^{\ast })\leq \alpha \frac{H_{\delta }(B_{n},S\left(
B_{n}\right) )[1+H_{\delta }(B^{\ast },T\left( B^{\ast }\right) )]}{%
1+H_{\delta }(B^{\ast },B_{n})}+H_{\delta }(T(B_{2n}),B_{n}),
\end{equation*}%
and on taking limit as $n\rightarrow \infty $ implies%
\begin{eqnarray*}
\lim \limits_{n\rightarrow \infty }H_{\delta }(B_{n},B^{\ast }) &\leq
&\lim \limits_{n\rightarrow \infty }\{ \alpha \frac{H_{\delta }(B_{n},S\left(
B_{n}\right) )[1+H_{\delta }(B^{\ast },T\left( B^{\ast }\right) )]}{%
1+H_{\delta }(B^{\ast },B_{n})}+H_{\delta }(T(B_{2n}),B_{n})\}, \\
\lim \limits_{n\rightarrow \infty }H_{\delta }(B_{n},B^{\ast }) &=&0,
\end{eqnarray*}%
which implies that $\lim \limits_{n\rightarrow \infty }B_{n}=B^{\ast }.$

\item[(vii)] If $M_{T,S}(B^{\ast },B_{n})=\dfrac{H_{\delta }(B_{n},T\left(
B^{\ast }\right) )[1+H_{\delta }(B_{n},T\left( B^{\ast }\right) )]}{%
1+H_{\delta }(B^{\ast },B_{n})}$ then%
\begin{equation*}
H_{\delta }(B_{n},B^{\ast })\leq \alpha H_{\delta }(B_{n},B^{\ast
})+H_{\delta }(T(B_{2n}),B_{n}),
\end{equation*}%
that is,%
\begin{equation*}
H_{\delta }(B_{n},B^{\ast })\leq \frac{1}{1-\alpha }H_{\delta
}(T(B_{2n}),B_{n}).
\end{equation*}%
On taking limit $n\rightarrow \infty $ implies%
\begin{eqnarray*}
\lim \limits_{n\rightarrow \infty }H_{\delta }(B_{n},B^{\ast }) &\leq &\frac{%
1}{1-\alpha }\lim \limits_{n\rightarrow \infty }H_{\delta }(T(B_{2n}),B_{n})
\\
&=&0,
\end{eqnarray*}%
which implies that $\lim \limits_{n\rightarrow \infty }B_{n}=B^{\ast }.$

\item[(viii)] If $M_{T,S}(B^{\ast },B_{n})=\dfrac{H_{\delta }(B_{n},T\left(
B^{\ast }\right) )[1+H_{\delta }(B^{\ast },T\left( B^{\ast }\right) )]}{%
1+H_{\delta }(B^{\ast },B_{n})}$ then%
\begin{equation*}
H_{\delta }(B_{n},B^{\ast })\leq \alpha \left( \frac{H_{\delta
}(B_{n},B^{\ast })}{1+H_{\delta }(B^{\ast },B_{n})}\right) +H_{\delta
}(T(B_{2n}),B_{n}),
\end{equation*}%
On taking limit $n\rightarrow \infty $ implies%
\begin{eqnarray*}
\lim \limits_{n\rightarrow \infty }H_{\delta }(B_{n},B^{\ast }) &\leq &\alpha
\lim \limits_{n\rightarrow \infty }\left( \frac{H_{\delta }(B_{n},B^{\ast })}{%
1+H_{\delta }(B^{\ast },B_{n})}\right) +\lim \limits_{n\rightarrow \infty
}H_{\delta }(T(B_{2n}),B_{n}), \\
&=&0,
\end{eqnarray*}%
which implies that $\lim \limits_{n\rightarrow \infty }B_{n}=B^{\ast }.$
Which completes the proof. $\square $\medskip 
\end{description}

\end{document}